\documentclass[12pt,reqno]{amsart}
\usepackage[margin=1in]{geometry}
\usepackage{amsmath,amssymb,amsthm}
\usepackage{mathrsfs}
\usepackage{xcolor}
\usepackage{hyperref}
\allowdisplaybreaks
\usepackage{quotes}

\numberwithin{equation}{section}
\newtheorem{Thm}{Theorem}[section]
\newtheorem{Def}{Definition}[section]
\newtheorem{Lem}{Lemma}[section]
\newtheorem{Pro}{Proposition}[section]

\newtheorem{Rem}{Remark}[section]

\def\mR{\mathbb{R}}
 
\def\mN{\mathbb{N}}

\def\mC{\mathbb{C}}

\def\R{\mathrm{R}}

\def\mcF{\mathbb {\mathcal F}}

\def\macB{\mathcal{B}}
\def\macS{\mathcal{S}}

\def\i0t{\int_0^t}

\let\originalleft\left
\let\originalright\right
\renewcommand{\left}{\mathopen{}\mathclose\bgroup\originalleft}
\renewcommand{\right}{\aftergroup\egroup\originalright}

\begin{document}
	
\newcommand{\Addresses}{{
		\bigskip
		\footnote{
	\noindent  \textsuperscript{1} U. S. Air Force Research Laboratory, Wright Patterson Air Force Base, Ohio 45433, U. S. A.
	
	\par\nopagebreak
	\noindent  \textsuperscript{2} NRC-Senior Research Fellow, National Academies of Science, Engineering and Medicine,  U. S. Air Force Research Laboratory, Wright Patterson Air Force Base, Ohio 45433, U. S. A.		
		\par\nopagebreak \noindent
		\textit{e-mail:} \texttt{provostsritharan@gmail.com} 	$^*$Corresponding author.
	\par\nopagebreak
\noindent  \textit{e-mail:} \texttt{saba.mudaliar@us.af.mil}

}}}

\title[Stochastic Quantization of Laser Propagation Models]{Stochastic Quantization of Laser Propagation Models \Addresses}
\author[ S. S. Sritharan and Saba Mudaliar]
{  Sivaguru S. Sritharan\textsuperscript{1,2*} and Saba Mudaliar\textsuperscript{1}}

\maketitle

\begin{abstract}
This paper identifies certain interesting mathematical problems of stochastic quantization type in the modeling of Laser propagation through turbulent media. In some of the typical physical contexts the problem reduces to stochastic Schr\"odinger equation with space-time white noise of Gaussian, Poisson and L\'evy type. We identify their mathematical resolution via stochastic quantization. Nonlinear phenomena such as Kerr effect can be modeled by stochastic nonlinear Schrodinger equation in the focusing case with space-time white noise. A treatment of stochastic transport equation, the Korteweg-de Vries Equation as well as a number of other nonlinear wave equations with space-time white noise is also given. Main technique is the S-transform (we will actually use closely related Hermite transform) which converts the stochastic partial differential equation with space time white noise to a deterministic partial differential equation defined on the Hida-Kondratiev white noise distribution space. We then utlize the inverse S-transform/Hermite transform known as the characterization theorem combined with the infinite dimensional implicit function theorem for analytic maps to establish local existence and uniqueness theorems for path-wise solutions of these class of problems. The particular focus of this paper on singular white noise distributions is motivated by practical situations where the refractive index fluctuations in propagation medium in space and time are intense due to turbulence, ionospheric plasma turbulence, marine-layer fluctuations, etc. Since a large class of partial differential equations that arise in nonlinear wave propagation have polynomial type nonlinearities, white noise distribution theory is an effective tool in studying these problems subject to different types of white noises.
\end{abstract}



\tableofcontents

\section{Introduction}
 Stochastic partial differential equations of the type 
\begin{equation}
	i\frac{\partial}{\partial t}\psi (x,t, \omega) + [\Delta +\Gamma(x,t,\omega)] \psi (x,t,\omega)  + F(\psi(x,t,\omega))=0,
\end{equation} 
where $\Gamma(\cdot,\cdot,\cdot)$ is a random field describing the fluctuations in the medium are often encountered in electromagnetic and acoustic wave propagation problems in random nonlinear media.  They are usually called either "paraxial equation" or "parabolic equation model" in the engineering literature \cite{Sp2002, Tap1977, Gustafsson2019}. As we define later, here $(\Omega, \mcF, m)$ is a complete probability space and $\omega$ is an $\macS'(\mR^{d})$-valued random variable. In this paper we will study Gaussian, Poisson and Levy type white noise models for  $\Gamma(\cdot,\cdot,\cdot)$. 
In the Gaussian case for example, $\Gamma(x,t,\omega)$ can also formally represented using generalized derivative of $\mR^{n+1}$-dimensional Brownian sheet $B(x,t,\omega)$:
\begin{equation}
	\Gamma(x,t,\omega)=\frac{\partial^{n+1}}{\partial x_{1},\cdots \partial x_{n}\partial t}B(x,t,\omega), \hspace{.1in} x=(x_{1},\cdots, x_{n}),
\end{equation}
where the Brownian sheet $B(x,t,\omega)$ has covariance:
\begin{equation}
	\langle B(x,t,\cdot), B(y,\tau,\cdot)\rangle =\Pi_{k=1}^{n}( x_{k}\wedge y_{k})( t\wedge \tau), \hspace{.1in} x=(x_{1},\cdots, x_{n}), y=(y_{1},\cdots, y_{n}).	
\end{equation}
The space-time Gaussian white noise has covariance formally: 
\begin{equation}
	\langle \Gamma(x,t,\cdot), \Gamma(y,\tau,\cdot)\rangle = \delta(x-y)\delta(t-\tau).
\end{equation}

 All three type noise structures will be developed in detail in later sections. In this paper we will provide a systematic treatment of Laser propagation in random media highlighting some of the well-known physical phenomena such as deep turbulence. In the simplest yet mathematically nontrivial case, the problem reduces to stochastic Schrodinger equation with space-time white noise and we discuss its relationship to stochastic quantization. Nonlinear phenomena such as Kerr effect can be modeled by stochastic nonlinear Schrodinger equation \cite{Sp2002} in the focusing case again with space-time white noise. Similar equations also arise in nonlinear fiber optics \cite{Agrawal2011}. We will also discuss a number of other stochastic partial differential equations  studied in the context of random wave propagation phenomena such as stochastic transport equation \cite{Ogawa1973, Funaki1979} and stochastic Korteweg-de Vries equation \cite{Bouard1998}. The quantization method we use is in the spirit of Wick expansions used in quantum field theory (see for example \cite{Simon1974}) and the white noise calculus (\cite{Hida1993, Kuo1996}) method for stochastic partial differential equations is described in \cite{Holden2010} which also gives a comprehensive discussion of the literature in this subject.
 \section{Derivation of the paraxial equation from the Maxwell equations}
 In this section we will give a heuristic derivation of the paraxial equation starting from the Maxwell equation. The paraxial equation we derive is a very well-known model widely used in the literature on acoustic wave and electromagnetic wave propagation in random media \cite{Tatarski1961, Papa1973, Tap1977, Sp2002, Gustafsson2019}. Let $E, H, B, D$ denote respectively the electric field, magnetizing field, magnetic field and the displacement field.The Maxwell equations are:
 \begin{equation}
 	\nabla \times E = \frac{\partial B}{\partial t}
 \end{equation}
 \begin{equation}
 	\nabla \times H = \frac{\partial D}{\partial t}
 \end{equation}
 \begin{equation}
 	D=\epsilon E \mbox{ and } B=\mu H.
 \end{equation}
 Here $\epsilon$ is the permitivity and $\mu$ is the permiability of the medium. 
 Taking curl of the first equation and taking time derivative of the second equation and substituting in the first we get upon using the vector identity $\nabla \times \nabla \times E = -\Delta E + \nabla (\nabla \cdot E)$ 
 \begin{equation}
 	\Delta E - \mu \frac{\partial^{2}}{\partial t^{2}}	(\epsilon E)= \nabla(\nabla \cdot E)
 \end{equation}
 In the absence of free charge $\nabla \cdot D=0$ and hence $E\cdot \nabla \epsilon + \epsilon \nabla \cdot E=0$ and substituting we get
 \begin{equation}
 	\Delta E - \mu \frac{\partial^{2}}{\partial t^{2}}	(\epsilon E)= -\nabla (E\cdot \nabla (log \epsilon))
 \end{equation}
 We have noting $\epsilon(x,t)=\epsilon_{0}n^{2}(x,t)$ where $n$ is the refractive index and $c^{2}=1/(\mu \epsilon_{0})$ with $c$ the light speed, we arrive at
 \begin{equation}
 	\Delta E(x,t) - \frac{1}{c^{2}} \frac{\partial^{2}}{\partial t^{2}}	(n^{2}(x,t) E(x,t))= -2\nabla (E(x,t)\cdot\nabla log (n(x,t))).	
 \end{equation}
 We assume that the time scale of fluctuations in the medium is much slower than the light speed and invoke further simplifications based on this assumption. Thus neglecting the right hand side and also $n^{2}$ term out of time derivative we arrive at
 \begin{equation}
 	\Delta E(x,t) - \frac{n^{2}(x,t)}{c^{2}} \frac{\partial^{2}}{\partial t^{2}}	 E(x,t)=0.
 \end{equation}
 Substituting a plane wave solution $E(x_{1},x_{2},x_{3},t)=\psi(x_{1},x_{2},x_{3},t)\exp(ik x_{3}-i\omega t)$ and neglecting the back-scatter term  $\frac{\partial^{2}\psi(x_{1},x_{2},x_{3},t)}{\partial x_{3}^{2}}$ (using simple scaling argument, see for example \cite{Papa1973}) we arrive at the paraxial equation:
 \begin{equation}
 	2ik \frac{\partial \psi}{\partial x_{3}}+[\Delta_{\bot} +k^{2} \{\frac{n^{2}(x_{1},x_{2},x_{3},t)}{n_{0}^{2}}-1\}]\psi=0.
 \end{equation}
 where $\Delta_{\bot} $ is the two dimensional Laplacian in variables $x_{2},x_{3}$. 
 Renaming the time-like variable $x_{3}$ as $t$ and suppressing the actual time variable $t$ we end up with the {\em two dimensional} linear or nonlinear stochastic Schr\"odinger equation:
 
 \begin{equation}
 	i\frac{\partial}{\partial t}\psi (x,t, \omega) + [\Delta +V(x,t,\omega, \psi)] \psi(x,t,\omega)=0.	
 \end{equation} 
 In this paper $V$ will be modeled as a space-time Gaussian white noise\cite{Hida1993,Kuo1996}, Poisson white noise \cite{Ito1988, ItoKubo1988, Benth1998} or Levy type white noise\cite{Holden2010} and the paraxial equation will then be framed as a stochastic quantization problem.
 
\section{Mathematical background on white noise calculus}
\subsection{White noise theory: Gaussian, Poisson and L\'evy}
In this section we will recall some of the basic elements of white noise stochastic calculus \cite{Hida1980, Hida1993, Kuo1996, Benth1998, Holden2010} needed in this paper. We will start by defining the Schwartz space $\macS=\macS(\mR^{d})$ of rapidly decaying real-valued  $C^{\infty}$ functions on $\mR^{d}$. This space equipped with the family of seminorms:
\begin{equation}
\| f\|_{k,\alpha}:=\sup_{x\in \mR^{d}}\{(1+\vert x\vert^{k})\vert\partial^{\alpha} f(x)\vert\}
\end{equation}
is a Frechet space \cite{Reed1980}. Here $k$ is a non-negative integer, $\alpha=\{\alpha_{1},\cdots, \alpha_{d}\}$ is a multi-index of non-negative integers $\alpha_{i},$ $ i=1,\cdots, d$ and 
\begin{equation}
\partial^{\alpha}f(x)=\frac{\partial^{\vert \alpha\vert}}{\partial x_{1}^{\alpha_{1}}\cdots \partial x_{d}^{\alpha_{d}}}f(x), \mbox{  where  } \vert \alpha\vert : =\alpha_{1}+\cdots +\alpha_{d}.
\end{equation}
The dual space of $\macS(\mR^{d})$ denoted $\macS':= \macS'(\mR^{d}) $ equipped with the weak topology is the space of tempered distributions. Let $(\Omega, \mcF, m)$ is a complete probability space and $\omega$ is an $\macS'(\mR^{d})$-valued random variable. The law $\mu$ of this $\macS'(\R^{d})$-valued random variable $\omega$ is characterized next. We recall the Bochner-Milnos theorem \cite{Reed1980} for the existence of probability measures on the Borel sets $\macB(\macS')$:
\begin{Thm}
A necessary and sufficient condition for the existence of a probability measure $\mu$ on $\macB(\macS')$ and a functional $F$ on $\macS$ such that
\begin{equation}
\int_{\macS'} e^{i <\omega, \phi>}	\mu(d\omega)=F(\phi), \hspace{.1in} \forall \phi \in \macS
\end{equation} 
is that $F$ satisfies the following three conditions:
\begin{enumerate}
	\item $F(0)=1$,
	\item $F$ is positive definite: $\sum_{j,l=1}^{n}z_{j}\bar{z}_{l}F(\phi_{j}-\phi_{l})\geq 0, \hspace{.1in} \forall z_{k}\in \mC, \forall \phi_{k}\in \macS, k=1,\cdots,n,$
	\item $F$ is continuous in the Frechet topology of $\macS$.	
\end{enumerate}	
\end{Thm}
We will also utilize three particular cases:\\
(i) Gaussian measure $\mu_{G}$:
\begin{equation}
\int_{\macS'} e^{i <\omega, \phi>}	\mu_{G}(d\omega)= \exp\{-\frac{1}{2}\|\phi\|_{L^{2}(\mR^{d})}^{2}\}, \hspace{.1in} \forall \phi \in \macS,	
\end{equation}
(ii) Poisson measure $\mu_{P}$:
\begin{equation}
	\int_{\macS'} e^{i <\omega, \phi>}	\mu_{P}(d\omega)= \exp\{\int_{\mR^{d}} (e^{i\phi(x)}-1)dx \}, \hspace{.1in} \forall \phi \in \macS,	
\end{equation}

(iii) Pure jump Levy measure $\mu_{L}$:
\begin{equation}
	\int_{\macS'} e^{i <\omega, \phi>}	\mu_{L}(d\omega)= \exp\{\int_{\mR^{d}} \Psi(\phi(x))dx \}, \hspace{.1in} \forall \phi \in \macS \mbox{  and  }
\end{equation}
\begin{equation}
\Psi(\phi)= \int_{\R^{d}}(e^{i<\phi,z> }-1-i<\phi,z>)	\nu(dz).
\end{equation}
\begin{Def}
\begin{enumerate}
	\item The continuous version of the mapping $\mR^{d}\ni x=(x_{1},\cdots, x_{d})\rightarrow B_{x}(\omega)\in L^{2}(\mu_{G})$ defined by
	\begin{equation}
		B_{x}(\omega)=\langle \omega, \Xi(x_{1})\times \cdots \times \Xi(x_{d})\rangle
	\end{equation}
is called the $d$-parameter Brownian motion (Brownian sheet). Here $\Xi(\cdot)\in L^{2}(\mR)$ is defined as:
\begin{equation}
\Xi(s) =
\left\{
\begin{array}{ll}
	\chi_{(0,s]}  & \mbox{if } s \geq 0 \\
	-\chi_{(s,0]} & \mbox{if } s < 0
\end{array}
\right.
\end{equation}
where $\chi$ is the usual indicator function.
\item The right continuous interger-valued mapping $\R^{d}\ni x=(x_{1},\cdots, x_{d})\rightarrow P_{x}(\omega)\in L^{2}(\mu_{P})$ defined by
\begin{equation}
	P_{x}(\omega)=\langle \omega, \Xi(x_{1})\times \cdots \times \Xi(x_{d})\rangle
\end{equation}
is called the $d$-parameter Poisson process. The mapping $\R^{d}\ni x=(x_{1},\cdots, x_{d})\rightarrow P_{x}(\omega)-\Pi_{i=1}^{d}x_{i}\in L^{2}(\mu_{P})$ is called compensated Poisson process.
\end{enumerate}
\end{Def}
We will now describe Wiener-Ito expansions:
\begin{Def}
Each $f\in L^{2}(\mu_{G})$ has an expansion in multiple (Brownian) Wiener integrals:
\begin{equation}
	f(\omega)=\sum_{n=0}^{\infty}\int_{\mR^{nd}}f_{n}(x)dB^{\otimes n}_{x}(\omega),
		\end{equation}
	where $f_{n}\in\hat{L^{2}(\R^{nd})}$ are deterministic symmetrized functions in $nd$ variables and
	\begin{equation}
		\|f\|^{2}_{L^{2}(\mu_{G})}=\sum_{n=0}^{\infty}n!\|f_{n}\|^{2}_{L^{2}(\mR^{nd})}.
	\end{equation}
Similarly
 $g\in L^{2}(\mu_{P})$ has an expansion in multiple (Poisson) Wiener integrals:
\begin{equation}
	g(\omega)=\sum_{n=0}^{\infty}\int_{\mR^{nd}}f_{n}(x)d(P_{x}(\omega)-\Pi_{i=1}^{d}x_{i})^{\otimes n}_{x}(\omega),
\end{equation}
where $g_{n}\in\hat{L^{2}(\mR^{nd})}$ are deterministic symmetrized functions in $nd$ variables and
\begin{equation}
	\|g\|^{2}_{L^{2}(\mu_{P})}=\sum_{n=0}^{\infty}n!\|g_{n}\|^{2}_{L^{2}(\mR^{nd})}.
\end{equation}
\end{Def}
We will now recall equivalent expansions in terms of Hermite and Charlier polynomials. For $n=1, 2, \cdots$ let $\zeta_{n}(x)\in \macS(\mR) $ be the Hermite function of order $n$:
\begin{equation}
\zeta_{n}(x):=\pi^{-1/4}((n-1)!)^{-1/2}	e^{-x^{2}/2}h_{n-1}(\sqrt{2x}), \hspace{.1in} x\in \mR,
\end{equation}
where $h_{n}(x)$ is the $n$-th Hermite polynomial defined by
\begin{equation}
	h_{n}(x):= (-1)^{n}e^{x^{2}/2}\frac{d^{n}}{dx^{n}}(e^{-x^{2}/2}), \hspace{.1in} x\in \mR, n=0,1,2, \cdots.
\end{equation}
It is well known that \cite{Reed1980} the sequence $\{\zeta_{n}\}_{n=1}^{\infty}$ forms an orthonomal basis for $L^{2}(\mR)$. Hence the family $\{\eta_{\alpha}\}$ of tensor products
\begin{equation}
\eta_{\alpha}=e_{\alpha_{1}\alpha_{2}\cdots \alpha_{d}}:=\zeta_{\alpha_{1}}\otimes \cdots \otimes 	\zeta_{\alpha_{d}}, \hspace{.1in} \alpha \in \mN^{d}
\end{equation}
forms an orthonomal basis for $L^{2}(\mR^{d})$. let us now assume that the family of all multi-indices $\alpha=(\alpha_{1}, \cdots, \alpha_{d})$ is given a fixed ordering 
\begin{equation}
(	\alpha^{(1)}, \alpha^{(2)}, \cdots \alpha^{(n)}, \cdots ),
\end{equation}
where $\alpha^{(k)}=(\alpha^{(k)}_{1}, \cdots, \alpha^{(k)}_{d})$ and denote $\eta_{k}=\eta_{\alpha^{(k)}}$. For a multi-index $\alpha=(\alpha_{1},\cdots, \alpha_{n})$ and $n\in \mN$ define the Hermite polynomial functionals as
\begin{equation}
H_{\alpha}(\omega):=\Pi_{j=1}^{n}h_{\alpha_{j}}(\langle \omega, \eta_{j}\rangle),
\end{equation}
and the Charlier polynomial functionals as
\begin{equation}
C_{\alpha}(\omega):=C_{\vert \alpha\vert}(\omega; \overbrace{\eta_{1},\cdots,\eta_{1}}^{\alpha_{1} \mbox{ times }},\cdots, \overbrace{\eta_{n},\cdots,\eta_{n}}^{\alpha_{n}\mbox{ times }}),
\end{equation}
with the convention
\begin{equation}
	C_{n}(\omega; e_{1}, \cdots, e_{n}):=\frac{\partial^{n}}{\partial \theta_{1}\cdots \partial \theta_{n}}e^{\{ \langle \omega, \ln{(1+\sum_{j=1}^{n}\theta_{j}e_{j})}\rangle -\sum_{j=1}^{n}\theta_{j}\int_{\mR^{d}}e_{j}(y)dy    \} }\vert_{\theta_{1}=\cdots=\theta_{n}=0}
\end{equation}
We also have the following equalities:
\begin{equation}
H_{\alpha} =\int_{\mR^{nd}} e^{\hat{\otimes}\vert \alpha\vert}d B_{x}^{\otimes \vert \alpha\vert }(\omega)
\end{equation}
and
\begin{equation}
	C_{\alpha} =\int_{\mR^{nd}} e^{\hat{\otimes}\vert \alpha\vert}d (P_{x}(\omega)-\Pi_{i=1}^{d}x_{i})_{x}^{\otimes \vert \alpha\vert }(\omega)
\end{equation}
Moreover, for any random functional $f(\omega)$ taking values on a separable Hilbert space $V$, and square integrable in $\mu_{G}$, $f\in L^{2}(\mu_{G};V)$ we have the Wiener Hermite polynomial chaos expansion \cite{Cameron1947}:
\begin{equation}
f(\omega)= \sum_{\alpha}a_{\alpha}H_{\alpha}(\omega), \hspace{.1in} a_{\alpha}\in V,
\end{equation}
with
\begin{equation}
	\|f\|^{2}_{L^{2}(\mu_{G};V)}=\sum_{\alpha}\alpha! \|a_{\alpha}\|_{V}^{2}, \mbox{  where } \alpha!=\alpha_{1}!\cdots \alpha_{n}!.
\end{equation}
Similarly, for any $g\in L^{2}(\mu_{P};V)$ we have the Wiener Charlier polynomial chaos expansion
\begin{equation}
	g(\omega)= \sum_{\alpha}b_{\alpha}C_{\alpha}(\omega), \hspace{.1in} b_{\alpha} \in V,
\end{equation}
with
\begin{equation}
	\|g\|^{2}_{L^{2}(\mu_{P};V)}=\sum_{\alpha}\alpha! \|b_{\alpha}\|_{V}^{2}, \mbox{  where } \alpha!=\alpha_{1}!\cdots \alpha_{n}!
\end{equation}
Note also that the correspondence between Gaussian and Poisson spaces ${\mathcal U}: L^{2}(\mu_{G};V)\rightarrow L^{2}(\mu_{P};V)$:
\begin{equation}
	{\mathcal U}\{\sum_{\alpha}a_{\alpha}H_{\alpha}(\omega)\} =\sum_{\alpha}a_{\alpha}C_{\alpha}(\omega),
\end{equation}
is unitary \cite{ItoKubo1988, Benth1997}. 
\subsection{Hida-Kondratiev spaces}
We will start with the characterization of the classical embedding of the Schwartz space in the space of tempered distributions as $\macS(\R^{d})\subset L^{2}(\mR^{d})\subset \macS'(\mR^{d})$ using Hermite Fourier coefficients \cite{Simon1971} and then define the Hida-Kondratiev spaces in an analogous way using Gaussian, Poisson and Levy measures. 
\begin{Thm}
(i) Let $\phi \in L^{2}(\mR^{d})$ so that we have the Fourier-Hermite expansion 
\begin{equation}
	\phi = \sum_{j=1}^{\infty} a_{j}\eta_{j},\hspace{.1in} \mbox{ where } a_{j}=(\phi, \eta_{j}), \hspace{ .1in} j=1,2, \cdots, 
\end{equation}
with 
\begin{equation}
	\eta_{j}:=\zeta_{\delta^{(j)}}=\zeta_{\delta_{1}^{(j)}}\otimes \cdots \zeta_{\delta_{d}^{(j)}}, \hspace{.1in} j=1,2, \cdots.
\end{equation}

Here $a_{j}$ are the Fourier coefficients of $\phi$ with respect to the tensor product Hermite functions $\eta_{j}$. Then $\phi \in S(\R^{d})$ if and only if 
\begin{equation}
	\sum_{j=1}^{\infty} a_{j}^{2}(\delta^{(j)})^{\gamma}<\infty,
\end{equation}
for all $d$-dimensional multi-indices $\gamma=(\gamma_{1},\cdots, \gamma_{d})$.

(ii) A distribution $T\in \macS'(\R^{d})$ is characterized by the expansion 
\begin{equation}
	T=\sum_{j=1}^{\infty} b_{j}\eta_{j}, \mbox{  with }
\end{equation}
\begin{equation}
	\sum_{j=1}^{\infty} b_{j}^{2}(\delta^{(j)})^{-\theta}<\infty,
\end{equation}
for some $d$-dimensional multi-index $\theta = (\theta_{1},\cdots, \theta_{d})$.
\end{Thm}
\begin{Def} Let $0\leq \rho\leq 1$. We say $f(\omega)=\sum_{\alpha} a_{\alpha}H_{\alpha} \in L^{2}(\mu_{G};V)$ belongs to Gaussian Hida-Kondratiev stochastic test function space $(\macS)_{G}^{\rho}(V)$ if
	\begin{equation}
\|f\|_{\rho,k}^{2}=\sum_{\alpha}\|a_{\alpha}\|_{V}^{2}(\alpha!)^{1+\rho}(2\mN)^{\alpha k}<\infty \hspace{.1in} \forall k \in \mN,
	\end{equation}
where
\begin{equation}
(2\mN)^{\alpha }=\Pi_{j=1}^{k}(2j)^{\alpha_{j}}, \hspace{.1in} \alpha= (\alpha_{1}, \cdots, \alpha_{k}).
\end{equation}
We say $f(\omega)=\sum_{\alpha} b_{\alpha}H_{\alpha} (\omega)\in L^{2}(\mu_{G};V)$ belongs to Gaussian Hida-Kondratiev stochastic distribution space $(\macS)_{G}^{-\rho}(V)$ if
\begin{equation}
	\sum_{\alpha}\|b_{\alpha}\|_{V}^{2}(\alpha!)^{1-\rho}(2\mN)^{-\alpha q}<\infty \mbox{ for some } q \in \mN.
\end{equation}
\end{Def}
This establishes the embedding
\begin{displaymath}
(\macS)^{1}_{X}(V)\subset (\macS)^{\rho}_{X}(V)\subset (\macS)^{+0}_{X}(V)\subset L^{2}(\mu_{X};V)\subset (\macS)^{-0}_{X}(V)\subset (\macS)^{-\rho}_{X}(V)\subset (\macS)^{-1}_{X}(V), 
\end{displaymath}
 with  $ X=G, P$ or $ L$.

{\bf Remark:}
The unitary correspondence  ${\mathcal U}: L^{2}(\mu_{G};V)\rightarrow L^{2}(\mu_{P};V)$:
\begin{equation}
	{\mathcal U}\{\sum_{\alpha}a_{\alpha}H_{\alpha}(\omega)\} =\sum_{\alpha}a_{\alpha}C_{\alpha}(\omega),
\end{equation}
can be extended as a unitary map ${\mathcal U}: (\macS)^{\rho}_{G}(V)\rightarrow (\macS)^{\rho}_{P}(V),$  $ -1\leq \rho \leq 1 $ (\cite{ItoKubo1988, Benth1998}).
\begin{Def} We have the following definitions of white noise processes.
\begin{enumerate} 
	\item A $d$-parameter Gaussian white noise is defined by 
\begin{equation}
W_{x}(\omega)=\sum_{k=1}^{\infty}\eta_{k}(x)H_{\varepsilon^{(k)}}(\omega), \hspace{.1in} x\in \mR^{d},
\end{equation}
where $\varepsilon^{(k)}$ is the multi-index with $1$ in $k$-th entry and zero otherwise.
\item A $d$-parameter compensated Poissonian white noise is defined by 
\begin{equation}
	\dot{P}_{x}(\omega)-1=\sum_{k=1}^{\infty}\eta_{k}(x)C_{\varepsilon^{(k)}}(\omega),  \hspace{.1in} x\in \mR^{d}.
\end{equation}
\end{enumerate}
\end{Def}
\begin{Lem}
We have
\begin{enumerate}
	\item $W_{x}(\omega)=\frac{\partial^{d}}{\partial x_{1}\cdots \partial x_{d}}B_{x}(\omega) \in (\macS)^{-\rho}_{G}$, $\rho\in [0,1]$,

	\item $\dot{P}_{x}(\omega)-1=\frac{\partial^{d}}{\partial x_{1}\cdots \partial x_{d}}(P_{x}(\omega)-\Pi_{i=1}^{d}x_{i}) \in (\macS)^{-\rho}_{P}$, $\rho\in [0,1]$.
\end{enumerate}

\end{Lem}

\subsection{Wick products and properties}
The Wick product is defined as below:
\begin{Def} The Wick product $F\diamond G$ of two elements of $(\macS)^{-1}(\mR^{n})$ is defined by:
	\begin{equation}
		F=\sum_{\alpha}a_{\alpha}H_{\alpha}, G=\sum_{\alpha}b_{\alpha}H_{\alpha} \in (\macS)^{-1}, \mbox{ with } a_{\alpha},b_{\alpha}\in \mR^{n},
	\end{equation}

\begin{equation}
F\diamond G :=\sum_{\alpha,\beta} a_{\alpha}\cdot b_{\beta} H_{\alpha+\beta} \in (\macS)^{-1}(\mR^{n}).
\end{equation}
\end{Def}
The fact that Wick product gives a distribution $F\diamond G \in  (\macS)^{-1}(\mR^{n})$ follows from the estimate below (see \cite{Holden2010}). We have $F=\sum_{\alpha}a_{\alpha}H_{\alpha}, G=\sum_{\beta}b_{\beta}H_{\beta}\in (\macS)^{-1}(\mR^{n})$ means there exists $q_{1}\in \mN$ such that

\begin{displaymath}
\sum_{\alpha}\vert a_{\alpha}\vert^{2}	(2\mN)^{-q_{1}\alpha} <\infty \hspace{.1in} \mbox{ and } \sum_{\beta}\vert b_{\beta}\vert^{2}	(2\mN)^{-q_{1}\beta} <\infty.
\end{displaymath}
Now rewriting 
\begin{equation}
F\diamond G :=\sum_{\alpha,\beta} a_{\alpha}\cdot b_{\beta} H_{\alpha+\beta} = \sum_{\gamma}(\sum_{\alpha +\beta =\gamma}a_{\alpha}\cdot b_{\beta}) H_{\gamma},
\end{equation}
and then setting  $C_{\gamma} =\sum_{\alpha +\beta =\gamma} a_{\alpha}\cdot b_{\beta} $ and $q=q_{1}+k$ we have by Cauchy-Schwartz inequality:
\begin{displaymath}
	\sum_{\gamma}(2 \mN)^{-q \gamma}\vert c_{\gamma}\vert^{2}\leq \left (\sum_{\gamma}(2 \mN)^{-k \gamma} \right ) \left (\sum_{\alpha}\vert a_{\alpha}\vert^{2}	(2\mN)^{-q_{1}\alpha}\right )\left (\sum_{\beta}\vert b_{\beta}\vert^{2}	(2\mN)^{-q_{1}\beta}   \right) <\infty.
\end{displaymath}
The first term on the right is finite for $k>1$ due to a Lemma by Zhang\cite{Zhang1992} and the other two terms are finite by definition of the distributions $F,G\in (\macS)^{-1}(\mR^{n})$.

\subsection{S-transform, Hermite transform and characterization theorems}
\begin{Def} Let $F=\sum_{\alpha}b_{\alpha}H_{\alpha}\in (S)_{G}^{-1}(V)$. Then the Hermite transform of $F$ denoted ${\mathcal H}_{G}F$ is defined :
\begin{equation}
	{\mathcal H}_{G}F=\sum_{\alpha}b_{\alpha}z^{\alpha} \in V   \hspace{.1in} \mbox{  (when convergent)},
\end{equation}
where $z=(z_{1},z_{2}\cdots )\in \mC^{\mN}$, $z^{\alpha}=z^{\alpha_{1}}\cdots z^{\alpha_{k}}$ for $\alpha=(\alpha_{1},\cdots, \alpha_{k}).$
\end{Def}
Similar statements for the Poisson and L\'evy cases. Remark 2.7, \cite{Zhang1992} clarifies the convergence of this series. 

\begin{Lem}
If $F,G\in (\macS)_{X}^{-1}(V)$, with $X=G,P,L$ (Gaussian, Poisson and Levy respectively) then 
\begin{equation}
	{\mathcal H}_{X}(F\diamond G)(z)={\mathcal H}_{X}F(z)\cdot {\mathcal H}_{X}G(z),
\end{equation}	
for all $z$ such that ${\mathcal H}_{X}F(z)$ and $ {\mathcal H}_{X}G(z)$ exist (convergent).
\end{Lem}

\begin{Lem} Suppose $g(z_{1},z_{2},\cdots)$ is a bounded analytic function on ${\bf B}_{q}(\delta)$ for some $\delta>0$, $0<q<\infty$ where
	\begin{equation}
	{\bf B}_{q}(\delta):= \left \{z=(z_{1},z_{2},\cdots) \in \mC_{0}^{N}; \hspace{.1in} \sum_{\alpha\ne 0}\vert z^{\alpha}\vert^{2}(2\mN)^{\alpha q}< \delta^{2}.\right \}
	\end{equation}
then there exists $F\in (\macS)_{G}^{-1}(V)$ and $D\in (\macS)_{P}^{-1}(V)$ such that ${\mathcal H}_{G}F=g={\mathcal H}_{P}D.$
\end{Lem}
See also the full statement of characterization theorem for $(\macS)^{-1}$ for the Gaussian as well as Poisson and L\'evy cases in \cite{Holden2010} (Theorem 2.6.11, Theorem 5.4.19). In this paper we will utilize a vector-valued version (Hilbert or Banach space valued) of the Hida-Kondratiev spaces, Hermite transform and inverse Hermite transform (Characterization theorem) and these extensions are straightforward generalizations of what is stated above.

\begin{Rem} Many authors have introduced ${\mathcal S}$-transform of $(\macS)^{-1}(V)$ -valued distributions \cite{Hida1993, Kuo1996} which is closely related to the Hermite transform:
\begin{equation}
	{\mathcal S}\phi(\zeta):=\int_{\macS^{'}}\phi(\omega)\exp^{\diamond }\langle \omega,\zeta\rangle d\mu(\omega),
\end{equation}
where $\exp^{\diamond }\langle \omega,\zeta\rangle =\sum_{0}^{\infty}\frac{1}{n!}\langle \omega,\zeta\rangle^{\diamond n}.$
It can be shown to be related to the Hermite transform as
\begin{equation}
	{\mathcal H}\phi(z_{1},z_{2},\cdots, z_{k})=	{\mathcal S}\phi(z_{1}\eta_{1}+\cdots+z_{k}\eta_{k})\hspace{.1in} \mbox{ for all } z_{1},z_{2},\cdots, z_{k} \in \mC^{k},
\end{equation}
in a suitable neighborhood. We however find it more convenient to work with the Hermite transform as pointed out in \cite{Holden2010}.
\end{Rem}
We also note for later use that the Hermite transform of the above white noises result in:
\begin{enumerate} 
	\item Hermite transform of $d$-parameter Gaussian white noise is given by: 
	\begin{equation}
		{\mathcal H}[W_{x}](z)=\sum_{k=1}^{\infty}\eta_{k}(x)z^{\varepsilon_{k}}.
	\end{equation}
	\item Hermite transform of $d$-parameter compensated Poisson white noise is given by 
	\begin{equation}
		{\mathcal H}[	\dot{P}_{x}-1](z)=\sum_{k=1}^{\infty}\eta_{k}(x)z^{\varepsilon_{k}}.
	\end{equation}
\end{enumerate}
Similar mathematical theory for pure jump L\'evy measures are given in Chapter 5 of Holden \cite{Holden2010} which parallel the above developments in polynomial expansions, distributional spaces as well as Hermite transforms which will be utilized in our paper as well.

\section{Kato theory of quasilinear abstract evolution equations }
All the problems considered in this paper, after Wick-quantization followed by Hermit (or $\macS $ )-transform, brought to a general class of deterministic quasilinear evolutions (complex in general) parameterized by an infinite sequence of complex numbers $z_{1}, z_{2},\cdots, $. We will thus recall some general results developed by T. Kato \cite{Kato1975, Kato1976} on quasilinear evolutions on a Banach space $X$:
\begin{equation}
\frac{d u}{dt}+A(t,u)u=0, \hspace{.1in} 0\leq t\leq T,
\end{equation}
\begin{equation}
	u(0)=u_{0}.
\end{equation}
Here $A(t,u)$ is an unbounded linear operator that nonlinearly depends on $t$ and $u$. Kato also points out in \cite{Kato1976} that an inhomogeneous equation with a right hand side:
\begin{equation}
	\frac{d u}{dt}+A(t,u)u= f(t,u)
\end{equation}
can be recast as a homogeneous problem above by redefining the variables. 
 Kato's theory covers parabolic problems where the linearized operator generates an analytic semigroup as well as hyperbolic problems \cite{Kato1970, Kato1973} where the linearized operator generates a $C_{0}$-semigroup which is not analytic. The main idea is to construct a mild solution $w\in C([0,T];X)$ for each $u\in C([0,T];X)$ for the linearized problem:
 \begin{equation}
 	\frac{d w}{dt}+A(t,u)w=0, \hspace{.1in} 0\leq t\leq T,
 \end{equation}
 \begin{equation}
 	w(0)=u_{0}.
 \end{equation}
 This defines a correspondence $u \rightarrow w =F(u)$ in C([0,T];X). Then use fixed point theory to construct the solution of the original quasilinear problem. 
  Kato's quasilinear theory was extended to the stochastic case using infinite dimensional Ito calculus in\cite{Fernando2015, Mohan2017}. In this paper we will extend Kato's theory to white noise calculus realm to treat problems with singular noises. We briefly summarize Kato's theory below \cite{Kato1975, Kato1976, Pazy1983}.
\begin{Def} Let $B$ be a subset of a Banach space $X$ and for every $0\leq t\leq T$ and $b\in B$ let $A(t,b)$ be the infinitesimal generator of a $C_{0}$-semigroup $S_{t,b}(s),$ $s\geq 0$ on $X$. The family of operators $\{A(t,b)\}, (t,b)\in [0,T]\times B$, is {\bf stable} if there are constants $M\geq 1$, $\omega$ such that the resolvent set 
	\begin{equation}
		\rho(A(t,b))\supset ]\omega, \infty],  \hspace{.1in} (t,b)\in [0,T]\times B
	\end{equation} 
and
\begin{equation}
	\|\Pi_{j=1}^{k} (\lambda-A(t_{j}, b_{j}))^{-1}\|\leq M (\lambda-\omega)^{-k}\hspace{.1in} \mbox{  for  } \lambda > \omega
\end{equation}
for every finite sequence $0 \leq t_{i}\leq t_{2}\cdots \leq T$, $b_{j}\in B, 1\leq j\leq k$.
\end{Def}
Stable family of generators $\{A(t,b)\}, (t,b)\in [0,T]\times B$ has stability estimate \cite{Kato1975}
\begin{equation}
	\|\Pi_{j=1}^{k} S_{t_{j},b_{j}}(s_{j})\| \leq M  e^{\omega \sum_{j=1}^{k} s_{j}}\hspace{.1in} \mbox{  for  } s_{j}\geq 0
\end{equation}
for every finite sequence $0 \leq t_{i}\leq t_{2}\cdots \leq T$, $b_{j}\in B, 1\leq j\leq k$.

\begin{Thm} 
Let $u_{0}\in Y$ and let $B= \{ v\in Y; \|v-u_{0}\|_{Y}\leq r\}, r>0 $.	Let the stable family of $C_{0}$-semigroup generators $\{A(t,b)\}, (t,b)\in [0,T]\times B$ satisfy the assumptions $H_{1}-H_{4}$ of \cite{Pazy1983} (Section 6.4) then there is a $T', 0 < T'\leq T$ such that the initial value problem 
\begin{equation}
	\frac{d u}{dt}+A(t,u)u=0, \hspace{.1in} 0\leq t\leq T,
\end{equation}
\begin{equation}
	u(0)=u_{0},
\end{equation}
has a unique mild solution $u\in C([0,T'];X)$ with $u(t)\in B$ for $t\in [0,T']$.
\end{Thm}
Details of the hypothesis of this theorem are discussed in \cite{Kato1975, Kato1976} and also \cite{Pazy1983}. Kato formulated this abstract theory to cover a large class of evolution system of physics and mechanics including the ones studied in the subsequent sections.

\subsection{Analytic Mappings between Banach Spaces and Inverse Mapping Theorem}

\begin{Def}
	A map $\Phi: Z_{1}\rightarrow Z_{2}$ between Banach spaces $Z_{1},Z_{2}$ is called analytic in $K_{\delta} = \left \{ u\in Z_{1}; \|u\|_{Z_{1}} <\delta \right \}$ if it can be deveoped in to a series of the form:
	\begin{equation}
		\Phi(u)=\sum_{k=0}^{\infty} \Phi_{k}(u,u, \cdots,u),
	\end{equation}
where $\Phi_{k}(\cdot, \cdots, \cdot): Z_{1}^{\otimes k}\rightarrow Z_{2}$ are symmetric multi-linear operators that are bounded:
\begin{displaymath}
\|\Phi_{k}(\cdot)\|=\sup\left \{\|\Phi_{k}(u, \cdots, u)\|_{Z_{2}}; \|u\|_{Z_{1}}\leq 1\right \}< \infty,
\end{displaymath}
and the series converges in $Z_{2}$-norm for $u\in K_{\delta}$ in the following sense:
\begin{displaymath}
	\sum_{k=0}^{\infty}\|\Phi_{k}\| \rho^{k} <\infty, \hspace{.1in} \mbox{ for } 0<\rho\leq \delta.
\end{displaymath}

\end{Def}
We will now state the analytic inverse function theorem \cite{Bourbaki1967, Vishik1988, Vallent1988}:
\begin{Thm}
Suppose that $\Phi:Z_{1}\rightarrow Z_{2}$ is an analytic map in a neighborhood of the origin $0\in K_{\delta}\subset Z_{1}$ and that its Frech\'et derivative at the origin  $D\Phi(0)$ is an isomorphism from $Z_{1}\rightarrow Z_{2}$. Then locally $\Phi$ has a unique inverse operator which is analytic in the neighborhood of $\Phi(0)\in Z_{2}$.
\end{Thm}
See Vallent\cite{Vallent1988} for the analytic version of the implicit function theorem as well.

\section{ First-order PDE (transport-type)}
\subsection{First order random transport equation with temporal white noise}
Transport problems arise in many applications including radiative transfer \cite{Chandrasekhar1947, Chandrasekhar1958, Chandrasekhar1960}. 
In this section we will consider first order random transport equation with temporal white noise studied by S. Ogawa \cite{Ogawa1973}, T. Funaki\cite{Funaki1979}, and H. Kunita \cite{Kunita1990} and briefly summarize their results.  Let $(\Omega, \Sigma, m)$ be a complete probability space. Consider

		\begin{equation}
			\frac{\partial}{\partial t}u(x,t,\omega)+ n(x,t,\omega)\frac{\partial }{\partial x}u(x,t,\omega) =c(x,t,\omega)u(x,t,\omega) + d(x,t, \omega).
		\end{equation}	
Here the wave speed $c(x,t,\omega)$ and forcing $ d(x,t, \omega)$ can be deterministic or random. We will first discuss the one dimensional Ogawa-Funaki temporal white noise model with $C,D$ deterministic:
	\begin{displaymath}
		\frac{\partial}{\partial t}u(x,t,\omega)+ \left \{\frac{d}{dt}\beta(t,\omega)+ b(t,x) \right \}\frac{\partial }{\partial x}u(x,t,\omega)
\end{displaymath}
\begin{equation}
		=c(x,t)u(x,t,\omega) + d(x,t).
\end{equation}
Here $\beta(t)$ is the 1-D Brownian motion.

Funaki also considered the n-dimensional scalar random transport model and we first write with Smooth noise:
\begin{displaymath}
			\frac{\partial}{\partial t}u(x,t,\omega)+\sum_{i=1}^{n} n_{i}(x,t,\omega)\frac{\partial }{\partial x_{i}}u(x,t,\omega)
\end{displaymath}
\begin{equation}
			=c(x,t)u(x,t,\omega) + d(x,t).
\end{equation}

The Funaki's temporal white noise model takes the form:
\begin{displaymath}
		\frac{\partial}{\partial t}u(x,t,\omega)+\sum_{i=1}^{n} \left \{\sum_{j}^{n} a_{ij}(t,x)\frac{d}{dt}\beta_{j}(t,\omega)+ b_{i}(x,t)\right \}\frac{\partial }{\partial x_{i}}u(x,t,\omega)
\end{displaymath}
\begin{equation}
		=c(x,t)u(x,t,\omega) + d(x,t).
\end{equation}
Here $\beta_{j}(t), j=1,\cdots, n$ are independent 1-D Brownian motions.

	 Define Stratenovich differential equation
		\begin{equation}
			dX_{t}=a(t, X_{t})\circ dB_{t}+b(t, X_{t})dt, \hspace{.1in} t\in (r, T),
		\end{equation}	
	
	\begin{equation}
		X_{r}=x.
	\end{equation}	
	This is equivalent to  Ito differential equation of the form:
		\begin{equation}
			dX_{t}=a(t, X_{t}) dB_{t}+\tilde{b}(t, X_{t})dt, \mbox{ where }
		\end{equation}
		\begin{equation}
			\tilde{b}(x,t)=b(t,x)+ \frac{1}{2}(a'a)(t,x),
		\end{equation}
		\begin{equation}
			(a'a)(t,x)_{i}=\sum_{k,j}\frac{\partial a_{ij}}{\partial x_{k}}	a_{kj}.
		\end{equation}

	\begin{Lem}
		For each $r, t$ with $  0\leq r\leq t\leq T$, $X(r,t,\cdot)$ is a homeomorphism of $\R^{n}$ on to $\R^{n}$.
	\end{Lem}
	
	 Define time reversed process 
		\begin{equation}
			Y(r,t,y)=X^{-1}(r,t,\cdot)(y), \hspace{.1in} 0\leq r\leq t\leq T, y\in \R^{n}.
		\end{equation}
Consider:	
	
	\begin{displaymath}
		\frac{\partial}{\partial t}u(x,t,\omega)+\sum_{i=1}^{n} \left \{\sum_{j}^{n} a_{ij}(t,x)\frac{d}{dt}\beta_{j}(t,\omega)+ b_{i}(x,t)\right \}\frac{\partial }{\partial x_{i}}u(x,t,\omega)
	\end{displaymath}
	\begin{equation}
		=c(x,t)u(x,t,\omega) + d(x,t),
	\end{equation}	
	\begin{equation}
		u(x,0)=\phi(x), \hspace{.1in} x\in G \mbox{ and } u(x,t)=\psi(x,t), \hspace{.1in} (x,t)\in \partial G\times [0,T].
	\end{equation}
The probabilistic solution to the stochastic transport equation by Funaki is:
		\begin{displaymath}
			u(x,t, \omega)= \left \{ \phi(Y(0,t,x)) +\psi(Y(\sigma(t,x),t,x))\right \}\exp\left \{\int_{0}^{t}c(s, Y(s,t,x))ds\right\}
		\end{displaymath}
		\begin{equation}
			+\int_{0}^{t} d(s, Y(s,t,x))\exp\left \{\int_{s}^{t}c(r, Y(r,t,x))dr\right\}ds.
		\end{equation}	
	Here $\sigma(t,x)$ is the exit time for the domain $G$.

	\subsection{Stochastic transport model with space-time white noise}
	We will begin with the following simple model (similar model with temporal white noise was considered by Funaki \cite{Funaki1979}) with space-time white noise $\Gamma(x,t)$ as characteristic speed:  
	\begin{equation}
		\frac{\partial}{\partial t}u(x,t, \omega)+\Gamma(x,t,\omega) \frac{\partial }{\partial x}u(x,t,\omega) =0, \hspace{.1in} (x,t)\in \R\times [0,T],
	\end{equation}	
	
	\begin{equation}
		u(x,0,\omega)=\phi(x), \hspace{.1in} x\in \R.
	\end{equation}
We quantize this problem as:
	\begin{equation}
	\frac{\partial}{\partial t}u(x,t,\omega)+\Gamma(x,t,\omega)\diamond  \frac{\partial }{\partial x}u(x,t,\omega) =0, \hspace{.1in} (x,t)\in \R\times [0,T],
\end{equation}	

\begin{equation}
	u(x,0,\omega)=\phi(x), \hspace{.1in} x\in \R.
\end{equation}
We take Hermite transform to get	(denoting $[{\mathcal H}u]:=\tilde{u}$ )
\begin{equation}
	\frac{\partial}{\partial t}\tilde{u}(x,t,z)+ [{\mathcal H}\Gamma](x,t,z)\frac{\partial }{\partial x}\tilde{u}(x,t,z) =0,
\end{equation}	
\begin{equation}
	\tilde{u}(x,0, z)=\phi(x).
\end{equation}
We write down the equation of characteristic:
\begin{equation}
	\frac{d\tilde{u}}{ds}= 	\frac{\partial \tilde{u}(x,t,z)}{\partial t}\frac{dt}{ds}+ \frac{\partial \tilde{u}(x,t,z) }{\partial x}\frac{dx}{ds}
\end{equation}
and arrive at 
\begin{equation}
	\frac{d \tilde{u}}{ds}=0, \hspace{.1in} 	\frac{dt}{ds}=1, \mbox{  and } \frac{dx}{ds}= [{\mathcal H}\Gamma](x,t,z).
\end{equation}	
Solving we get $x=\zeta_{t}(x_{0})=x_{0}+ \int_{0}^{t}[{\mathcal H}\Gamma](x,r,z)dr$ as the equation of characteristics and $\tilde{u}(x,0,z)=\phi(x_{0})=\phi (\zeta_{t}^{-1}(x))$: 
\begin{equation}
	\tilde{u}(x,t,z)=\phi(\zeta_{t}^{-1}(x))=\phi(x-\int_{0}^{t}[{\mathcal H}\Gamma](x,r,z)dr).
\end{equation}
Hence
\begin{equation}
	u(x,t,\omega)={\mathcal H}^{-1}\phi(\zeta_{t}^{-1}(x))={\mathcal H}^{-1}\left (\phi(x-\int_{0}^{t}[{\mathcal H}\Gamma](x,r,z)dr)\right ).
\end{equation}
We will now turn to the multidimensional stochastic transport model with space-time white noise:
\begin{displaymath}
	\frac{\partial}{\partial t}u(x,t,\omega)+\sum_{i=1}^{n} \left \{\sum_{j}^{n} a_{ij}(t,x)\Gamma_{j}(x,t,\omega)+ b_{i}(x,t)\right \}\frac{\partial }{\partial x_{i}}u(x,t,\omega)
\end{displaymath}
\begin{equation}
	=c(x,t)u(x,t,\omega) + d(x,t).
\end{equation}	 
We will also study the quantized random transport model:
\begin{displaymath}
	\frac{\partial}{\partial t}u(x,t,\omega)+\sum_{i=1}^{n} \left \{\sum_{j}^{n} a_{ij}(t,x) \Gamma_{j}(x,t,\omega)\right \}\diamond \frac{\partial }{\partial x_{i}}u(x,t,\omega) +\sum_{i}^{n} b_{i}(x,t)\frac{\partial }{\partial x_{i}}u(x,t,\omega) 
\end{displaymath}
\begin{equation}
	=c(x,t)u(x,t,\omega) + d(x,t).
\end{equation}
 Taking the Hermite transform: For $z\in \mC^{\mN}$, and denoting ${\mathcal H} u $ as $\tilde{u}$, we get
\begin{displaymath}
	\frac{\partial}{\partial t}\tilde{u}(x,t,z)+\sum_{i=1}^{n} \left \{\sum_{j}^{n} a_{ij}(t,x) [{\mathcal H}\Gamma]_{j}(x,t,z) +  b_{i}(x,t)\right \} \frac{\partial }{\partial x_{i}}\tilde{u}(x,t,z) 
\end{displaymath}
\begin{equation}
	=c(x,t)\tilde{u}(x,t,z) + d(x,t),
\end{equation}
This is a deterministic linear hyperbolic equation (for a fixed $z$ ) and we can solve it similar to the simple model presented above and the apply inverse Hermite transform to obtain the solution for the original stochastic transport equation with space-time white noise as characteristic speed coefficient and spatial white noise initial data. With the help of the deterministic results in \cite{Diperna1989} we can also obtain the solvability theorem for stochastic transport equation with white noise characteristic speeds and spatial white noise initial data as
\begin{Pro} Suppose that $a,b\in L^{1}(0,T;L^{1}_{\mbox{loc}}(\mR^{n}))$, $c,d\in L^{1}(0,T;L^{\infty}(\mR^{n}))$ and the initial data is a spatial white noise with $u_{0}\in (\macS)^{-1}(L^{\infty}(\mR^{n}))$, and the stochastic characteristic speed coefficients $\Gamma_{j}\in (\macS)^{-1}(L^{\infty}(\mR^{n}\times [0,T]))$, $ j=1,\cdots,n$. Then there exists a unique solution $u\in (\macS)^{-1}(L^{\infty}(0,T;L^{\infty}(\mR^{n}))$.
\end{Pro}

We now consider a multidimensional quasilinear stochastic transport model with space-time white noise:
\begin{displaymath}
	\frac{\partial}{\partial t}u(x,t,\omega)+\sum_{i=1}^{n} \left \{\sum_{j}^{n} a_{ij}(u( t,x, \omega))\Gamma_{j}(x,t,\omega)+ b_{i}(u(x,t, \omega))\right \}\frac{\partial }{\partial x_{i}}u(x,t,\omega)
\end{displaymath}
\begin{equation}
	=c(u(x,t,\omega)),
\end{equation}	 
where $a_{ij}, b_{i}, c$ are all polynomials in $u$.  We quantize this problem as:
\begin{displaymath}
	\frac{\partial}{\partial t}u(x,t,\omega)+\sum_{i=1}^{n} \left \{\sum_{j}^{n} (a_{ij}(u( t,x, \omega))^{\diamond})\diamond \Gamma_{j}(x,t,\omega)+ b_{i}(u(x,t, \omega))^{\diamond}\right \}\diamond \frac{\partial }{\partial x_{i}}u(x,t,\omega)
\end{displaymath}
\begin{equation}
	= c(u(x,t,\omega))^{\diamond},
\end{equation}	 
Applying Hermite transform we get
\begin{displaymath}
\frac{\partial}{\partial t}\tilde{u}(x,t,z)+\sum_{i=1}^{n} \left \{\sum_{j}^{n} a_{ij}(\tilde{u}( t,x,z))[{\mathcal H}\Gamma]_{j}(x,t,z)+ b_{i}(\tilde{u}(x,t, z))\right \} \frac{\partial }{\partial x_{i}}\tilde{u}(x,t,z)
\end{displaymath}
\begin{equation}
= c(\tilde{u}(x,t,z)).
\end{equation}
 Extending the deterministic result to the Hermite transformed problem may require some smoothing of the white noise term $\Gamma$ using an operator of the form $(I-\Delta_{x,t})^{-\gamma}$ for a suitable $\gamma \geq 0$.
The local solvability of this deterministic partial differential equation for a given fixed $z$ can be obtained from Kato's theory \cite{Kato1975, Kato1976} as $\tilde{u}\in C([0,T'];H^{s}(\mR^{n}))$, $s>n/2+1$ and from this we can use the characterization theorem for Hermite transform to deduce:
\begin{Pro}
For $u_{0}\in (\macS)^{-1}( H^{s}(\mR^{n}))$, $s>n/2+1$ and $\Gamma_{j}(\cdot,\cdot,\cdot) \in ({\macS})^{-1}(L^{\infty}_{\mbox{loc}}(\mR^{n}\times \mR^{+}))$,  $ j=1,\cdots,n$, there is a  unique local-in-time generalized solution for the Gaussian white noise forced quasilinear transport equation $u\in (\macS)^{-1}(C([0,T'];H^{s}(\mR^{n})))$. For the Poisson and L\'evy cases the unique correspondence map ${\mathcal U} $ gives a unique solution $u\in (\macS)^{-1}(L^{\infty}([0,T'];H^{s}(\mR^{n})))$.
\end{Pro} 
 
 Global in time generalized solution (that accommodates shock waves) can be built by applying Kruzkov theory \cite{Kruzkov1970} to the Hermite transformed problem with small bounded variation norm initial data:. 
\begin{Pro}
	For $u_{0}\in (\macS)^{-1}\left ( L^{\infty}_{\mbox{loc}}(\mR^{n})\cap BV(\mR^{n})\right )$, with sufficiently small norm, and $\Gamma(\cdot,\cdot,\cdot) \in ({\macS})^{-1}(L^{\infty}_{\mbox{loc}}(\mR^{n}\times \mR^{+}))$ there is a  unique global-in-time generalized solution for the Gaussian/Poisson/L\'evy white noise forced quasilinear transport equation $u\in (\macS)^{-1}(L^{\infty}([0,T];L^{\infty}_{\mbox{loc}}(\mR^{n})))$. 
\end{Pro} 
\section{Stochastic nonlinear wave equations}

\subsection{Stochastic Korteweg De Vries equation}
Let us consider the Korteweg De Vries  equation with white noise initial data:
\begin{equation}
	\frac{\partial}{\partial t}\varphi (x,t, \omega)+ \varphi (x,t, \omega)\frac{\partial}{\partial x}\varphi (x,t, \omega)+\frac{\partial^{3}}{\partial x^{3}}\varphi (x,t, \omega) =0,
\end{equation}
\begin{equation}
	\varphi(x,0,\omega)=\Gamma(x,\omega).
\end{equation}
 We quantize this problem as
\begin{equation}
	\frac{\partial}{\partial t}\varphi (x,t, \omega)+ \varphi (x,t, \omega)\diamond \frac{\partial}{\partial x}\varphi (x,t, \omega)+\frac{\partial^{3}}{\partial x^{3}}\varphi (x,t, \omega) =0.
\end{equation}

Applying Hermite transform results in:
\begin{equation}
	\frac{\partial}{\partial t}\tilde{\varphi} (x,t, z)+ \tilde{\varphi} (x,t, z)\frac{\partial}{\partial x}\tilde{\varphi }(x,t, z)+\frac{\partial^{3}}{\partial x^{3}}\tilde{\varphi} (x,t, z) =0.
\end{equation}
\begin{equation}
\tilde{	\varphi}(x,0,z)=[{\mathcal H} \Gamma](x,z).
\end{equation}
Note that the Hermite transformed KDV equation will have infinite number of conservation laws as in P. D. Lax theory for deterministic KDV equation \cite{Lax1968} some of which would be finite depending on the smoothness of noise . However, for the case of singular white noise, values of these conservation laws will be all infinity.  

Tsutsumi \cite{Tsutsumi1983} has studied deterministic Korteweg-De Vries equation with bounded Radon measures as initial data (see also Bourgain \cite{Bourgain1997} for the space periodic case) and proven existence of a generalized solution and we can use this theorem the above Hermite transformed problem for a fixed $z$ followed by inverse Hermite transform to conclude that:
\begin{Pro}
For the Gaussian/Poisson/L\'evy white noise initial data $\varphi(\cdot,0,\cdot)=\Gamma(\cdot,\cdot)\in (\macS )^{-1}(L^{\infty}_{\mbox{loc}}(\mR))$, there exists a unique solution $u$ to the stochastic KDV equation such that:
\begin{displaymath}
u\in (\macS)^{-1}\left ( L^{\infty}(0,\infty;H^{-1}(\mR))\cap L^{2}((0,T)\times (-R,R))\right ),\hspace{.1in} \mbox{ for any } T, R>0.	
\end{displaymath}
\end{Pro}

\subsection{Stochastic Benjamin-Ono equation}
The stochastic Benjamin-Ono equation is modeled as:
\begin{equation}
	\frac{\partial}{\partial t}\varphi (x,t, \omega)+ \varphi (x,t, \omega)\frac{\partial}{\partial x}\varphi (x,t, \omega)+{\bf H}[\frac{\partial^{2}}{\partial x^{2}}\varphi (x,t, \omega)] =0,
\end{equation}
where ${\bf H}[\cdot]$ is the Hilbert transform defined as
\begin{equation}
	{\bf H}[f](x):=\mbox{ PV }\frac{1}{\pi} \int_{-\infty}^{\infty}\frac{f(y)}	{x-y}dy = \hspace{.1in} {\mathcal F}^{-1}\left (i \cdot (
\mbox{ sign}\zeta ){\mathcal F}(f)(\zeta)\right ),
\end{equation}
where ${\mathcal F}$ 
 and  ${\mathcal F}^{-1}$ denote Fourier transform and its inverse respectively.  This problem is supplied with random initial data:
 \begin{equation}
 \varphi(x,0)=\Gamma(x, \omega).	
 \end{equation}
 
 We quantize this equation as
 \begin{equation}
 	\frac{\partial}{\partial t}\varphi (x,t, \omega)+ \varphi (x,t, \omega)\diamond \frac{\partial}{\partial x}\varphi (x,t, \omega)+{\bf H}[\frac{\partial^{2}}{\partial x^{2}}\varphi (x,t, \omega)] =0,
 \end{equation}
Applying Hermite transform we get:
 \begin{equation}
	\frac{\partial}{\partial t}\tilde{\varphi} (x,t, z)+ \tilde{\varphi} (x,t, z)\frac{\partial}{\partial x}\tilde{\varphi }(x,t, z)+{\bf H}[\frac{\partial^{2}}{\partial x^{2}}\tilde{\varphi} (x,t, z)] =0,
\end{equation}
\begin{equation}
	\tilde{\varphi}(x,0,z)={\mathcal H}[\Gamma](x,z)].
\end{equation}
We can utilize the sharp results of T. Tao \cite{Tao2004} along with a suitable smoothing of the noise by an operator of the form $(I-\Delta)^{-\gamma}$ to conclude that there is a unique solution to the above problem for fized $z$ as $\tilde{\varphi}\in C([0,T];H^{s})$ for $s\geq 1$. Hence using inverse Hermite transform we obtain:

\begin{Pro}
For the Gaussian white noise initial data $\varphi(\cdot,0,\cdot)=\Gamma(\cdot,\cdot)\in (\macS )^{-1}(H^{s}(\mR))$, $s\geq 1$,  there exists a  unique solution to stochastic Benjamin-Ono equation as $\varphi \in (\macS)^{-1}(C([0,T];H^{s})),$ $ s\geq 1.$ For Poisson and L\'evy cases the unique correspondence map ${\mathcal U}$ gives unique solution $\varphi \in (\macS)^{-1}(L^{\infty}([0,T];H^{s})),$ $ s\geq 1.$
\end{Pro}

\section{Stochastic reaction diffusion equation and quantization}
Parisi and Wu \cite{Parisi1981} initiated the subject of stochastic quantization \cite{Damgaard1987} which takes the view that (Euclidean) quantum fields can be constructed by studying stochastic partial differential equations. Stochastic reaction diffusion equation with space-time Gaussian noise has been studied by many authors \cite{Faris1982} with results such as construction of invariant measures \cite{Benzi1989, Gatarek1996} as well as pathwise strong solutions\cite{Albeverio2001, DaPrato2003}. In this section we will treat this class of problems with white noise theory.

\subsection{Stochastic heat equation with multiplicative noise and the KPZ Equation}

Kardar, Parisi and Zhang \cite{Kardar1986} studied the stochastic model (which has come to be known as the KPZ equation):
\begin{equation}
\frac{\partial h}{\partial t}=\nu \Delta h + \frac{\lambda}{2}(\nabla h)^{2} + \Gamma(x,t).
\end{equation}
As pointed out in this paper, in the case of smooth noise the above stochastic PDE can be formally transformed in to two other well-known stochastic PDEs (ignoring some constants). In fact $v=-\nabla h$ gives the stochastic Burgers equation:
\begin{equation}
	\frac{\partial v}{\partial t} +\lambda v\cdot \nabla v =\nu \Delta v +\nabla \Gamma(x,t),
\end{equation}
and the transform $\varphi =\exp ((\frac{\lambda}{2\nu})h)$ converts the KPZ equation to

\begin{equation}
\frac{\partial}{\partial t}\varphi (x,t, \omega) =\nu \Delta \varphi (x,t,\omega) + \frac{\lambda}{2\nu }\varphi (x,t,\omega)\Gamma(x,t,\omega),
\end{equation}
with spatial white noise initial data:
\begin{displaymath}
\varphi (x,0, \omega)=\Gamma_{0}(x,\omega).	
\end{displaymath}
We quantize this problem as (after setting $\nu=\sigma^{2}/2$ and $\lambda=\sigma^{2}$ for simplicity):
\begin{equation}
	\frac{\partial}{\partial t}\varphi (x,t, \omega) =\frac{1}{2}\sigma^{2}\Delta \varphi (x,t,\omega) + \varphi (x,t,\omega)\diamond \Gamma(x,t,\omega).
\end{equation}
Applying the Hermite transform gives
\begin{equation}
	\frac{\partial}{\partial t}\tilde{\varphi} (x,t, z) -[\frac{1}{2}\sigma^{2}\Delta  + [{\mathcal H}\Gamma](x,t,z)]\tilde{\varphi} (x,t,z)=0,
\end{equation}
\begin{displaymath}
	\tilde{\varphi} (x,0, z)=[{\mathcal H}\Gamma_{0}](x,z).
\end{displaymath}
We write the solution formally using the propagator $ {\mathcal Z}(t,r)$ generated by \cite{ Pazy1983} the unbounded time dependent operator $ \frac{1}{2}\sigma^{2}\Delta +[{\mathcal H}\Gamma](x,t,z) $:

\begin{equation}
	\tilde{\varphi} (x,t, z)= {\mathcal Z}(t,0)	[{\mathcal H}\Gamma_{0}](x,z).
\end{equation}
The solution is also probabilistically expressed by the Feynman-Kac formula using a Brownian motion $B_{t}$ independent of $\Gamma(x,t,\omega)$ and $\Gamma_{0}(x,\omega)$:
\begin{equation}
\tilde{\varphi} (x,t, z)	=E_{x}\left [ [{\mathcal H}\Gamma_{0}](\sigma B_{t},z) \exp \left \{\int_{0}^{t} [{\mathcal H}\Gamma](\sigma B_{s},t-s,z)ds \right \}\right ].
 \end{equation}
Alternatively, we can consider the solution as a fixed point problem for the heat semigroup:
\begin{equation}
	\tilde{\varphi} (x,t, z)= e^{\frac{1}{2}\sigma^{2}t\Delta}[{\mathcal H}\Gamma_{0}](x,z)	+ \int_{0}^{t}e^{\frac{1}{2}\sigma^{2}(t-\tau)\Delta}{\mathcal H}[\Gamma](x,\tau,z)\tilde{\varphi} (x,\tau, z)d\tau.
\end{equation}
Heat equation with singular and non-autonomous potentials have been studied in the literature \cite{Gulisashvili2004, Gulisashvili2006} and for a fixed $z$ as necessary using smoothing operator of the form $(I-\Delta_{x,t})^{-\gamma}$ for the noise we can conclude that the above Hermite transformed problem has a unique solution in $\tilde{\varphi}\in C([0,T];L^{\infty})$ and hence by inverse Hermite transform we obtain:
\begin{Pro}
For the Gaussian white noise $\Gamma(\cdot,\cdot,\cdot)\in (\macS )^{-1}(L^{\infty}_{\mbox{loc}}(\mR^{d}\times \mR))$, and initial data $\varphi(\cdot,0)\in L^{\infty}_{\mbox{loc}}(\mR^{d})$,  there exists a unique solution to the stochastic heat equation $ 		\varphi \in (\macS)^{-1}(C([0,T];L^{\infty}(\mR^{d}))).$ For the Poisson and L\'evy white noises the unique correspondence map ${\mathcal U}$ gives  a unique solution to the stochastic heat equation $ 		\varphi \in (\macS)^{-1}(L^{\infty}([0,T];L^{\infty}(\mR^{d}))).$
 	
 The solution to the quantized problem is also probabilistically expressed by the Feynman-Kac formula:
 \begin{equation}
 \varphi(x,t, \omega)	=E_{x}\left [ \Gamma_{0}(\sigma B_{t},\omega)\diamond \exp^{\diamond} \left \{\int_{0}^{t} \Gamma(\sigma B_{s},t-s,\omega)ds \right \}\right ],	
 \end{equation}
where the Wick exponential of $X\in (\macS)^{-1}$ is defined by:
\begin{displaymath}
	\exp^{\diamond}X =\sum_{0}^{\infty}\frac{1}{n!}X^{\diamond n}.
\end{displaymath}
 \end{Pro}

We also note here that instead of the multiplicative noise term if we have a stochastic heat equation with space-time white noise as an additive noise term or initial data as a spatial white noise then after Hermite transform we end up with heat equation with singular initial data or forcing term and the problem can be resolved using the fundamental solution of the heat equation followed by inverse Hermite transform.

\subsection{Stochastic nonlinear heat equation with white noise Initial Data}

We will now consider:
\begin{equation}
	\frac{\partial}{\partial t}\varphi (x,t, \omega)+ \varphi (x,t, \omega)^{p}= \Delta \varphi (x,t,\omega),\hspace{.1in} p=2,3, x\in \mR^{n},
\end{equation}
\begin{displaymath}
\varphi (x,0, \omega)=\Gamma(x,\omega),  x\in \mR^{n}.
\end{displaymath}
We will quantize this equation as:
\begin{equation}
	\frac{\partial}{\partial t}\varphi (x,t, \omega)+\varphi (x,t, \omega)^{\diamond p} = \Delta \varphi (x,t,\omega),
\end{equation}
Applying the Hermite transform gives
\begin{equation}
	\frac{\partial}{\partial t}\tilde{\varphi} (x,t, z) -[\Delta -\tilde{\varphi} (x,t,z)^{p-1}]\tilde{\varphi} (x,t,z) =0,
\end{equation}
\begin{equation}
\tilde{\varphi}(x,0,z)=[{\mathcal H}\Gamma](x,z).
\end{equation}
Nonlinear heat equations with measure initial data has been studied in the literature with positive as well as negative results. When the initial data is a Dirac measure it has been shown in \cite{Brezis1981} that this problem with the nonlinearity selected above $(p=2$ and $3)$ has no solution in any space dimension $n\geq 1$. This means that we need to smooth the measure by an operator of the form $(I-\Delta)^{-\gamma}$ so that we can use results for slightly less singular but still measure data such as that presented in \cite{Brezis1996} with initial data in $L^{q}, 1\leq q <\infty$. This provides a unique short time solution to the Hermite transformed problem as $\tilde{\varphi}\in C([0,T'];L^{q})$ and we then inverse Hermite transform to obtain: 
\begin{Pro}
There exists a unique solution $\varphi \in (\macS)^{-1}(C([0,T'];L^{q}))$ to the stochastic nonlinear heat equation with spatial white noise initial data $\Gamma(\cdot,\cdot) \in  (\macS)^{-1}(L^{q})$.
 
\end{Pro}

\section{Stochastic linear and nonlinear Schr\"odinger equations with space-time white noise}
In this section we will address the stochastic linear and nonlinear Schr\"odinger equations with space-time white noise that arise in the Laser propagation problems as discussed earlier. Once again we will utilize the Hermite transform to convert the problem to deterministic linear and nonlinear Schrodinger equations parameterized by an infinite sequence of complex variables $z$. There is a wealth of literature on linear Schr\"odinger semigroup with a range of potentials \cite{Simon1982, Berezin1991} and nonlinear Schr\"odinger equations \cite{Bourgain1999, Sulem1999} which will enable the solvability of the deterministic problems obtained by Hermite transforms as we discuss below.

\subsection{Strichartz estimates}
We recall the following estimates for the Schr\"odinger free propagator \cite{Strichartz1977,Keel1998}:
\begin{Lem}
	For $(q,r)$ and $(\tilde{q},\tilde{r})$ such that for $d\geq 1$ both exponent pair satisfying:
	\begin{displaymath}
	2\leq q,r \leq \infty, \frac{1}{q}=\frac{d}{2} (\frac{1}{2}-\frac{1}{r})\hspace{.1in} (q,r,d)\neq (2,\infty,2),
	\end{displaymath}
we have:
\begin{equation}
	\| e^{it\Delta}u_{0}\|_{L^{q}_{t}L^{r}_{x}(\mR\times \mR^{d})} \lesssim \|u_{0}\|_{L^{2}_{x}(\mR^{d})},
\end{equation}
and
\begin{equation}
\|\int_{0}^{t}e^{i(t-s)\Delta}F(s,\cdot)ds\|_{L^{q}_{t}L^{r}_{x}(\mR\times \mR^{d})}\lesssim \|F\|_{L^{\tilde{q}'}_{t}L^{\tilde{r}'}_{x}(\mR\times \mR^{d})}.
\end{equation}
\end{Lem}

\subsection{Stochastic linear Schr\"odinger equation with additive space-time white noise}
\begin{equation}
	i\frac{\partial}{\partial t}\psi (x,t, \omega) + \Delta \psi (x,t,\omega) = \Gamma(x,t,\omega), \hspace{.1in} (x,t)\in \mR^{d}\times \mR^{+},
\end{equation}
with spatial white noise initial data:
\begin{equation}
\psi (x,0, \omega)=	\Gamma_{0}(x,\omega),\hspace{.1in} x\in \mR^{d}.
\end{equation}
Applying Hermite transform we get the Schr\"odinger equation in free space:
\begin{equation}
	i\frac{\partial}{\partial t}\tilde{\psi} (x,t, z) + \Delta \tilde{\psi} (x,t,z) = {\mathcal H}[\Gamma](x,t,z), \hspace{.1in} x\in \mR^{d},  z\in \mC^{\mN},
\end{equation}
\begin{equation}
	\tilde{\psi} (x,0, z)=	{\mathcal H}[\Gamma_{0}](x,z),\hspace{.1in} x\in \mR^{d}, z\in \mC^{\mN}.
\end{equation}
The solution is written formally using the free Schr\"odinger propagator $e^{it\Delta}$ as
	\begin{equation}
\tilde{\psi} (x,t, z)= e^{it\Delta}	{\mathcal H}[\Gamma_{0}](x,z)	-i \int_{0}^{t}e^{i(t-\tau)\Delta}{\mathcal H}[\Gamma](x,\tau,z)d\tau.
	\end{equation}

The d-dimensional Schr\"odinger kernel is:
\begin{equation}
	K_{t}(x)=\frac{1}{(4\pi it)^{d/2}}e^{i \frac{ \vert x\vert^{2}}{4t}}.
\end{equation}
Hence we also have:
\begin{equation}
\tilde{	\psi}(x,t)	=\frac{1}{(4\pi it)^{d/2}}\left [\int_{\R^{d}}	e^{i \frac{ \vert x-y\vert^{2}}{4t}}{\mathcal H}[\Gamma_{0}](y,z)dy-i \int_{0}^{t}\int_{\R^{d}}e^{i \frac{ \vert x-y\vert^{2}}{4(t-\tau)}}{\mathcal H}[\Gamma](x,\tau,z)dyd\tau \right ].
\end{equation}
This equation makes sense for a fixed $z\in \mC^{\mN}$ if ${\mathcal H}[\Gamma_{0}](x,z)\in L^{2}(\mR^{d})$ and  ${\mathcal H}[\Gamma](\cdot,\cdot,z)\in L^{q}_{t}L^{r}_{x}(\mR\times \mR^{d})$ due to the Strichartz estimates recalled above. Hence inverse Hermite transform gives:
\begin{Pro}
For initial data noise distribution $\Gamma_{0}\in (\macS)^{-1}(L^{2}(\mR^{d}))$ and noise forcing $\Gamma\in (\macS)^{-1}(L^{q}_{t}L^{r}_{x}(\mR\times \mR^{d}))$ there exists a unique solution $\psi\in (\macS)^{-1}(L^{q}_{t}L^{r}_{x}(\mR\times \mR^{d}))$ to the stochastic Schr\"odinger equation.  
\end{Pro}

\subsection{Stochastic linear Schr\"odinger equation with multiplicative space-time white noise}

 Consider the linear stochastic Schr\"odinger equation with multiplicative space-time noise: For $\omega \in \Omega, x\in \R^{d},d\geq 1, t>0$
\begin{equation}
	i\frac{\partial}{\partial t}\psi (x,t, \omega) + [ \Delta +\Gamma(x,t,\omega)]\psi(x,t,\omega)=0,
\end{equation}
and spatial white noise initial data:
\begin{equation}
	\psi (x,0, \omega)=	\Gamma_{0}(x,\omega),\hspace{.1in} x\in \mR^{d}.
\end{equation}
Here the time-like coordinate $t$ is the propagation direction, $\varphi$ is the (complex) electric field and the potential $V$ in general depends on the refractive index of the medium. 

We consider the quantized linear Schr\"odinger equation with multiplicative white noise:
\begin{equation}
	i\frac{\partial}{\partial t}\psi (x,t, \omega) + \Delta \psi (x,t,\omega) +  \psi(x,t,\omega)\diamond \Gamma(x,t,\omega)=0.
\end{equation}
Applying Hermite transform we get the Schr\"odinger equation with a potential $V(x,t,z)={\mathcal H} [\Gamma](x,t,z)$:

\begin{equation}
	i\frac{\partial}{\partial t}\tilde{\psi }(x,t, z) + \left ( \Delta  +  {\mathcal H} [\Gamma](x,t,z)\right )\tilde{\psi}(x,t,z)=0.
\end{equation}
\begin{equation}
	\tilde{\psi} (x,0, z)=	{\mathcal H}[\Gamma_{0}](x,z),\hspace{.1in} x\in \mR^{d}, z\in \mC^{\mN}.
\end{equation}
We write the solution formally using the propagator $ {\mathcal Z}(t,r)$ generated by \cite{Simon1982, Pazy1983} the unbounded Hamiltonian time dependent $ (\Delta + V(x,t,z))$:

	\begin{equation}
	\tilde{\psi} (x,t, z)= {\mathcal Z}(t,0)	\tilde{\psi} (x,0, z).
\end{equation}
Alternatively, we can consider the solution as a fixed point problem for the free propagator:
	\begin{equation}
	\tilde{\psi} (x,t, z)= e^{it\Delta}\tilde{\psi} (x,0, z)	-i \int_{0}^{t}e^{i(t-\tau)\Delta}{\mathcal H}[\Gamma](x,\tau,z)\tilde{\psi} (x,\tau, z)d\tau.
\end{equation}
Schr\"odinger equations with time dependent unbounded singular potential have been studied in the literature \cite{Yajima2011} where it is shown that the two parameter propagator $ {\mathcal Z}(t,r)$ is unitary in $L^{2}(\mR^{2})$ (actually K. Yajima's results holds in any dimension.) With an introduction to a suitable smoothing for the noise in the form of $(I-\Delta_{x,t})^{-\gamma}$ we can fit the Hermite transformed problem above in Yajima's framework and deduce a unique solution $\tilde{\psi} \in C([0,T];L^{2}(\mR^{2}))$ and hence inverse Hermite transform gives:

\begin{Pro}
For Gaussian initial data noise distribution $\Gamma_{0}\in (\macS)^{-1}(L^{2}(\mR^{d}))$ and Gaussian white noise forcing $\Gamma\in (\macS)^{-1}(L^{\infty}_{\mbox{loc}}(\mR\times \mR^{d}))$, there exists a unique solution to the stochastic Schr\"odinger equation with multiplicative space-time white noise as $\psi\in (\macS)^{-1}(C([0,T];L^{2}(\mR^{2})))$. For  Gaussian and L\'evy cases unique correspondence map ${\mathcal U}$ gives a unique solution $\psi\in (\macS)^{-1}(L^{\infty}([0,T];L^{2}(\mR^{2})))$
\end{Pro}

\subsection{Nonlinear Schr\"odinger equation with multiplicative space-time white noise}
Nonlinear defocusing cubic Stochastic Schr\"odinger equation with multiplicative space-time white noise:
\begin{equation}
	i\frac{\partial}{\partial t}\psi (x,t, \omega) + \Delta \psi (x,t,\omega)+ \left (  \Gamma(x,t,\omega)-\vert\psi\vert^{2}\right )\psi(x,t,\omega)=0, \hspace{.1in} x\in \mR^{3}, t>0
\end{equation}
and spatial white noise initial data:
\begin{equation}
	\psi (x,0, \omega)=	\Gamma_{0}(x,\omega),\hspace{.1in} x\in \mR^{3}.
\end{equation}
We consider the quantized nonlinear Schrodinger equation with multiplicative white noise:
\begin{equation}
	i\frac{\partial}{\partial t}\psi (x,t, \omega) + \Delta \psi (x,t,\omega)+ \left( \Gamma(x,t,\omega) - \psi(x,t,\omega)\diamond \psi^{*}(x,t,\omega)\right)\diamond \psi(x,t,\omega)=0.
\end{equation}
Here (in the Gaussian white noise case) we take
\begin{displaymath}
	\psi (x,t, \omega)=\sum_{\alpha}\Phi_{\alpha}(x,t)H_{\alpha}(\omega) \mbox{  and its complex conjugate } 	\psi^{*} (x,t, \omega)=\sum_{\alpha}\Phi^{*}_{\alpha}(x,t)H_{\alpha}(\omega),
\end{displaymath}
Applying Hermite transform we get a deterministic nonlinear Schr\"odinger equation:
\begin{equation}
i	\frac{\partial}{\partial t}\tilde{\psi} (x,t, z) +[(\Delta +{\mathcal H}[\Gamma](x,t,z)  -  \tilde{ \psi}(x,t,z)\tilde{ (\psi^{*})}(x,t,z)]\tilde{ \psi}(x,t,z)=0,
\end{equation}
where
\begin{displaymath}
\tilde{\psi} (x,t, \omega)=\sum_{\alpha}\Phi_{\alpha}(x,t)z^{\alpha} \mbox{  and  } \tilde{ (\psi^{*})}(x,t,z)	=\sum_{\alpha}\Phi^{*}_{\alpha}(x,t)z^{\alpha}.	
\end{displaymath}
\begin{equation}
	\tilde{\psi} (x,0, z)=	{\mathcal H}[\Gamma_{0}](x,z),\hspace{.1in} x\in \mR^{3}, z\in \mC^{\mN}.
\end{equation}
Nonlinear Schr\"odinger equation with time dependent potential of the above type (with fixed $z$ ) is studied in \cite{Carles2011}. With an introduction to a suitable smoothing for the noise in the form of $(I-\Delta_{x,t})^{-\gamma}$ the Hermite transformed problem above in the framework of \cite{Carles2011} (in particular Assumption 1.3 regarding the potential in that paper) and deduce a unique solution $\tilde{\psi} \in C([0,T];L^{2}(\mR^{3}))\cap L^{8/3}([0,T];L^{4}(\mR^{3}))$ and hence inverse Hermite transform gives:

\begin{Pro}
For Gaussian initial data noise distribution $\Gamma_{0}\in (\macS)^{-1}(L^{2}(\mR^{d}))$ and Gaussian noise forcing $\Gamma\in (\macS)^{-1}(L^{\infty}_{\mbox{loc}}(\mR\times \mR^{d}))$, there exists a unique solution to the stochastic nonlinear Schr\"odinger equation with multiplicative space-time white noise as $\psi\in (\macS)^{-1}(C([0,T];L^{2}(\mR^{3}))\cap L^{8/3}([0,T];L^{4}(\mR^{3})))$. For Poisson and L\'evy cases unique correspondence map ${\mathcal U}$ gives a unique solution $\psi\in (\macS)^{-1}(L^{\infty}([0,T];L^{2}(\mR^{3}))\cap L^{8/3}([0,T];L^{4}(\mR^{3})))$
\end{Pro}

\section{Concluding remarks}
In this paper we have casted a number of Laser propagation problems in random media in the white noise distribution theory framework to stimulate further research in their mathematical structure. Some of the most prominent problems are selected for our discussion and a number of other nonlinear wave equations that arise in Laser interaction with plasma such as the Schrodinger-Hartree equation \cite{Hartree1957}, Zakharov system \cite{Zakharov1972} and Davey-Stewartson equation \cite{Davey1974, Carbonaro2010} can also be treated by the method initiated here when dealing with stochastic medium effects. We briefly indicate below how one may proceed with Wick quantization in these well-known nonlinear wave problems and for simplicity we formulate them with spatially random initial data. \\

(1) Random (quantized) Schr\"odinger-Hartree equation:
\begin{equation}
	i\frac{\partial}{\partial t}\varphi + \Delta \varphi = \pm [\vert x\vert^{n}\star (\varphi^{*}\diamond \varphi )]\diamond \varphi,  \hspace{.1in} x\in \mR^{d}, \hspace{.1in} 0< n <d,
\end{equation}
\begin{equation}
	\varphi(x,\omega,0)=\Gamma(x,\omega).
\end{equation}
Upon Hermite transform this system (for a fixed $z$) will result in the usual deterministic Schr\"odinger-Hartree system measure initial data. We can apply a smoothing to the noise in the form $(I-\Delta)^{-\gamma}$ and then use a solvability theorem such as in \cite{Hyakuna2019} and utilize the analytic implicit function theorem \cite{Bourbaki1967, Vallent1988} and  inverse Hermite transform to deduce the solvability:
\begin{Pro}
For $0<\gamma <\mbox{ min }(2,n)$, and Gaussian/Poisson/L\'evy white noise initial data $\Gamma(\cdot,\cdot)\in (\macS)^{-1}(L^{2}(\mR^{n}))$, there is a unique solution to the (quantized) Schr\"odinger-Hartree equation $\varphi \in (\macS)^{-1}\left (C(\mR;L^{2}(\mR^{n}))\cap L^{8/\gamma}_{\mbox{loc}}(\mR;L^{\frac{4n}{2n-\gamma}}(\mR^{n}))\right)$. 
\end{Pro}

\vspace{.1in}

(2) Random (quantized) Zakharov system in $\mR^{d+1}$, $d=2,3$:
\begin{equation}
	i\frac{\partial}{\partial t}\varphi + \Delta \varphi = \varphi \diamond n,
\end{equation}
\begin{equation}
	[\frac{\partial^{2}}{\partial t^{2}}-\Delta ]n = -\Delta (\varphi^{*}\diamond \varphi),
\end{equation}
\begin{equation}
	\varphi(x,\omega,0)=\Gamma_{1}(x,\omega), n(x,\omega,0)=\Gamma_{2}(x,\omega), \mbox{ and  } \frac{\partial}{\partial t} n(x,\omega,0)=\Gamma_{3}(x,\omega).
\end{equation}
Upon Hermite transform this system (for a fixed $z$) will result in the usual deterministic Zakharov system measure initial data. We may have to apply a smoothing to the noise in the form $(I-\Delta)^{-\gamma}$ before we can start with a suitable solvability theorem such as in \cite{Ginibre1997,Chen2021} and utilize the analytic implicit function theorem \cite{Bourbaki1967, Vallent1988} and  inverse Hermite transform \cite{Holden2010} to deduce the solvability of the stochastic Zakharov model.
\begin{Pro} Suppose the initial data Gaussian white noise distribution satisfy $(\Gamma_{1},\Gamma_{2},\Gamma_{3})\in (\macS)^{-1}\left (H^{1/2}(\mR^{d}))\times L^{2}(\mR^{d})\times H^{-1}(\mR^{d})\right )$ then there exists a unique local-in-time solution to stochastic Zakharov equation with white noise initial data such that 
	\begin{equation}
	(\varphi, n, \partial_{t} \varphi)\in (\macS)^{-1}\left ( C([0,T];H^{1/2}(\mR^{d})\times L^{2}(\mR^{d})\times H^{-1}(\mR^{d}))\right ).
\end{equation}
\end{Pro}

\vspace{.1in} 

(3) Random (quantized) Davey-Stewartson system in $\mR^{2+1}$:
The system was derived by Davey and Stewartson \cite{Davey1974} and see \cite{Ghidaglia1990} for a rigorous study. We consider the Wick-quantized problem:
\begin{equation}
i\frac{\partial}{\partial t}u + \delta\frac{\partial^{2}}{\partial x^{2}}u+\frac{\partial^{2}}{\partial y^{2}}u=\chi (u^{*}\diamond u)\diamond\varphi + b u\diamond \frac{\partial}{\partial x}\varphi,
\end{equation}
\begin{equation}
	\frac{\partial^{2}}{\partial x^{2}}\varphi +m\frac{\partial^{2}}{\partial y^{2}}\varphi =\frac{\partial}{\partial x}(u^{*}\diamond u),
\end{equation}
\begin{equation}
	u(x,y,\omega,0)=\Gamma(x,y,\omega).
\end{equation}
The four parameters $\delta, \chi, b, m$ are real, $\vert \delta \vert =\vert \chi \vert =1$. The system is classified as elliptic-elliptic, elliptic-hyperbolic, hyperbolic-elliptic and hyperbolic-hyperbolic according to the respective signs of $(\delta, m): (+.+), (+,-), (-,+), (-,-)$. Upon Hermite transform this system (for a fixed $z$) will result in the usual deterministic Davey-Stewartson system measure initial data and as discussed in the paper we can start with a suitable solvability theorem such as in \cite{Ghidaglia1990, Villamizar2012} and utilize the analytic implicit function theorem \cite{Bourbaki1967, Vallent1988} and  inverse Hermite transform \cite{Holden2010} to deduce the solvability of the stochastic model: 
\begin{Pro}
For the Gaussian/Poisson/L\'evy initial data white noise distribution $\Gamma\in (\macS)^{-1} (L^{2}(\mR^{2}))$, there exists a unique local in time solution in the elliptic-elliptic and hyperbolic-elliptic cases ($m>0)$:
\begin{equation}
u\in (\macS)^{-1}\left ( C([0,T^{*}[; L^{2}(\mR^{2}))\cap L^{4}((0,T^{*})\times \mR^{2})\right ),
\end{equation}
\begin{equation}
\nabla \varphi  \in (\macS)^{-1}\left (L^{2}((0,T^{*})\times \mR^{2})\right ).
\end{equation}
\end{Pro}

\vspace{.3in}
{\bf Acknowledgment}: The first author's research has been supported by the U. S. Air Force Research Laboratory through the National Research Council Senior Research Fellowship of the National Academies of Science, Engineering and Medicine.


\end{document}